\documentclass{article}
\usepackage[utf8]{inputenc}
\usepackage[margin = 1in]{geometry}

\usepackage{amsmath}
\usepackage{bm}
\usepackage{graphicx}
\usepackage{xcolor}
\usepackage{hyperref}
\usepackage{dsfont}
\def\b{\ensuremath\mathbf}

\usepackage[toc,page]{appendix}

\usepackage[numbers]{natbib}

\usepackage{graphicx}% Include figure files
\usepackage{dcolumn}% Align table columns on decimal point
\usepackage{bm}% bold math
\usepackage[utf8]{inputenc}
\usepackage[T1]{fontenc}
\usepackage{mathptmx}
\usepackage{etoolbox}

\usepackage{xcolor}
\usepackage{ dsfont }

\usepackage{hyperref}

\usepackage{amsfonts}

\def\b{\ensuremath\mathbf}
\title{Chaotic heteroclinic networks as  models of switching behavior in biological systems}

\author{Megan Morrison$^1$ and Lai-Sang Young$^1$\\[.1in]
$^1$ Courant Institute of Mathematical Sciences, New York University, New York, NY 10012
}

\begin{document}

\maketitle

\begin{abstract}
Key features of biological activity can often be captured by transitions between a finite number of
semi-stable states that correspond to behaviors or decisions. We present here a broad class of
dynamical systems that are ideal for modeling such activity. The models we propose are 
chaotic heteroclinic networks with nontrivial intersections of stable and unstable manifolds. 
Due to the sensitive dependence on initial conditions, transitions between states are 
seemingly random. Dwell times, exit distributions, and 
other transition statistics can be built into the model through geometric design and can be controlled
by tunable parameters. To test our model's ability to simulate realistic biological phenomena, we
turned to one of the most studied organisms, {\it C. elegans}, well known for its limited
behavioral states. We reconstructed experimental data 
from two laboratories, demonstrating the model's ability to quantitatively reproduce dwell times
and transition statistics under a variety of conditions. Stochastic switching between dominant
states in complex dynamical systems has been extensively studied and is often modeled 
as Markov chains.
As an alternative, we propose here a new paradigm, namely chaotic heteroclinic networks 
generated by deterministic rules (without the necessity for noise). Chaotic heteroclinic
 networks can be used to model systems with arbitrary architecture and size without
a commensurate increase in phase dimension. They are highly flexible and able to capture 
 a wide range of transition characteristics that can be adjusted through control parameters.
\end{abstract}

\maketitle

\section{\label{sec:level1}Introduction}

Biological systems are invariably highly complex, involving very large numbers of 
components
or substances interacting in complicated ways. Two examples are cortical and 
metabolic networks. Without major simplification, analytical studies of 
detailed biological models seem hopelessly out of reach. Not all biological 
models can be dimension reduced, but when a system has a finite number of 
dominant states that are semi-stable, phenomenological models focusing 
on transitions between these states can be constructed and analyzed. 
One way to reduce to this more tractable situation
is to sufficiently constrain the circumstances or scope of  study, thereby limiting
 the set of relevant behaviors. When they exist, dominant states may be identified 
 through empirical observation. They can also be deduced from methods of
 data-driven modeling such as PCA, HMM, Koopman approximation, or
 machine learning \cite{brunton_data-driven_2019, kutz_dynamic_2016, westhead_hidden_2017, alla_nonlinear_2017, kato_global_2015, taghia_uncovering_2018}.
 
This paper is about dynamical models of systems  with a finite number of 
semi-stable ``attractor states" around which the system stabilizes briefly before 
transitioning to another state. We are interested in situations where trajectories are 
not periodic but seemingly random; which state will occur next has the appearance 
of being unpredictable, as do transition times.
Switching dynamics of this kind have been observed in many biological 
settings (see Discussion). The attractor states in question typically represent fundamental
 behaviors or processes \cite{kato_global_2015,chang_endogenous_2021,wiltschko_mapping_2015, sharma_point_2018, rabinovich_transient_2008}, and the branching is often associated with decisions. 
Our motivating examples are primarily from biology, but we
 consider here general dynamical systems with the characteristics above.  

Continuous-time Markov chains (the states of which correspond to the attractor
states above) are the simplest natural  models 
that come to mind for describing switching dynamics, but the Markovian 
assumption is  strong. Biological events often have some degree of history dependence.
What an animal will do next, for example, may depend on what it did last and
how long it has been in its present state. It would be useful to have tractable
models that permit a broader set of dynamical behaviors. 
 
 We propose that 
{\it chaotic heteroclinic networks} are excellent phenomenological models of 
switching dynamics. Given a homoclinic loop or heteroclinic cycle, it has been
known since the days of Poincar\'e that transversal intersections of stable and
unstable manifolds give rise to complicated behavior\cite{palis_hyperbolicity_1995, poincare1890probleme, poincare1899methodes}. 
The qualitative picture is clarified by Smale, who showed that this complicated 
behavior included the existence of horseshoes \cite{smale_differentiable_1967}. 
Here we take these ideas one step further.
We demonstrate that by manipulating the eigenvalues of the saddle fixed points 
in a heteroclinic network together with the geometry of stable and unstable manifolds 
and their intersections, one can impose {\it quantitative control} on the switching dynamics, including
branching probabilities, dwell times, and other more detailed characteristics. 
These new tools have enabled us to design large classes of chaotic dynamical systems 
and to tailor transition statistics in heteroclinic networks in order to match observed behavior 
in biological modeling.

As proof of concept, we constructed chaotic heteroclinic networks 
to reproduce experimental {\it C. elegans} data. \textit{C. elegans} 
are a well-studied model organism in biology due to their simple anatomy, simple nervous system, few observed behaviors, and the relative ease of performing experimental measurements of behavioral
and neural activity \cite{kato_global_2015,nguyen_whole-brain_2016, prevedel_simultaneous_2014, schrodel_brain-wide_2013, white_structure_1986}. It has been shown that such activity 
can be represented in low-dimensional PCA space where semi-stable states correspond to stereotyped behaviors (or fictive behaviors) with characterizable probabilities of transition between these states \cite{kato_global_2015, nichols_global_2017, linderman_hierarchical_2019}. \textit{C. elegans} transition between different behaviors seemingly at random; their tendency to reside in various dynamical and behavioral states can be modified with experimental conditions, such as oxygen levels, genetic strain, and developmental stage \cite{nichols_global_2017}. The stochastic switching dynamics of \textit{C. elegans}, modulated by experimental conditions and corresponding to different locations in neural activity space, is an ideal system to demonstrate how heteroclinic networks can generate stochastic switching behavior observed in biological systems.

Researchers have modeled \textit{C. elegans}' activity with a variety of different paradigms. 
Markov models can accurately capture the observed switching statistics but do not shed light
on how transitions are generated \cite{roberts_stochastic_2016, gallagher_geometry_2013, linderman_hierarchical_2019}; it has also been observed that target states vary with 
dwell times, pointing to the need for more detailed modeling \cite{nichols_global_2017}. 
Previous dynamical systems models of \textit{C. elegans} activity used control inputs 
as the mechanism for inducing state changes \cite{fieseler_unsupervised_2020, morrison_nonlinear_2021-1},
but state changes in \textit{C. elegans} appear to occur mostly spontaneously
and not in direct response to external stimuli. 
The  model we propose has the capability to address all of these issues.

Seemingly random switching from one pattern of behavior to another
with no apparent trigger is a motif seen throughout biology and neuroscience.
We propose that chaotic heteroclinic networks with quantitative control of
transition properties may be suitable phenomenological models for 
this type of behavior.

\section{Heteroclinic networks: general theory and quantitative control}

This section contains the theory part of the paper.
We begin by reviewing some basic ideas connected with heteroclinic dynamics.
Then we specialize to 2D, discussing first the case 
of flows before perturbing their time-$t$-maps to produce branching behavior
in chaotic heteroclinic networks, the main objects of interest in this paper. 
This is followed by examples of quantitative control on the switching dynamics, 
techniques that will be used to build phenomenological models of {\it C. elegan} 
behavior.

\medskip \noindent
{\it Basic definitions}

\smallskip
Consider an ordinary differential equation (ODE) 
\begin{align}
    \dot{x} = f(x), \quad x \in \mathds{R}^d, \label{eq:dxdt=fx}
\end{align}
where $d \ge 2$ is an arbitrary integer and $f$ is assumed to be smooth. 
The flow associated with Eq.~\ref{eq:dxdt=fx} is denoted by $\varphi_t$, so that
if $x(t)$ is the solution of the ODE with initial condition $x(0) = x_0 \in \mathds{R}^d$,
then $\varphi_t(x_0):=x(t)$ is the trajectory of the flow starting from $x_0$.
Let $p$ be an equilibrium point of $\varphi_t$ i.e. $f(p)=0$. We say $p$ is a hyperbolic 
fixed point  if none of the eigenvalues $\lambda$ of $Df(p)$ lies on the imaginary axis. 
If the real parts of $\lambda$ are $<0$ (resp. $>0$) for all $\lambda$, then $p$ is a sink
(resp. a source); if some are $>0$ and some are $<0$, then $p$ is a saddle fixed point.
Through each hyperbolic fixed point of saddle type are stable and unstable manifolds, defined to be
\begin{align}
W^s(p) &= \{ x \in \mathds{R}^d: \varphi_t(x) \rightarrow p \ \ \text{as} \ \ t \rightarrow \infty\},\\
W^u(p) &= \{ x \in \mathds{R}^d: \varphi_t(x) \rightarrow p \ \ \text{as} \ \ t \rightarrow -\infty\}
\end{align}
respectively. A point $x$ is called a \textit{homoclinic point} of $p$ if $x \in W^u(p) \cap W^s(p)$.
If $p_1$ and $p_2$ are distinct hyperbolic fixed points of saddle type, then $x$ is a  
\textit{heteroclinic point} from $p_1$ to $p_2$ if $x \in W^u(p_1) \cap W^s(p_2)$, i.e.,
$\varphi_t(x) \rightarrow p_1$ as $t \to -\infty$ and $\varphi_t(x) \rightarrow p_2$ as $t \to \infty$.
A point $x$ is called a \textit{tranverse homoclinic point} of $p$ if $W^u(p)$ meets $W^s(p)$
transversally at $x$. Transverse heteroclinic points are defined similarly.

In discrete time, i.e., for dynamical systems generated by iterating a map $F: \mathds{R}^d
\to \mathds{R}^d$ which we assume to be smooth and invertible, the corresponding objects
are defined analogously: e.g. a fixed point $p \in \mathds{R}^d$ is called a hyperbolic fixed
point of saddle type if the spectrum $\Sigma$ of $DF(p)$ is such that $\Sigma \cap \{|z|=1\}
= \emptyset$ and $\Sigma \cap \{|z|<1\}, \ \Sigma \cap \{|z|>1\} \ne \emptyset$.
Stable and unstable manifolds 
are defined as above with $n \in \mathbb Z$ in the place of $t$, and 
a point $x \in  W^u(p_1) \cap W^s(p_2)$ is called a heteroclinic point 
from $p_1$ to $p_2$.

The existence of transverse homoclinic points is well known to lead to complicated dynamics,
an observation that goes back to Poincar\'e. The qualitative picture is developed by Smale, 
who showed that homoclinic points are accompanied by
the existence of horseshoes, a hallmark of chaos. For a more systematic treatment,
see Refs.~\cite{palis_hyperbolicity_1995, palis2012geometric, meiss2007differential, guckenheimer2013nonlinear}.

\medskip \noindent
{\it Heteroclinic networks for flows in 2D}

\smallskip
The continuous-time picture in $\mathbb R^2$ is simple
because for saddle fixed points, $W^u$ and
$W^s$ are 1D, and if $x$ is a heteroclinic point from $p_1$ to $p_2$, then the orbit of $x$,
$\mathcal O(x):=\{\varphi_t(x), t \in \mathbb R\}$, coincides with one of the branches of
$W^u(p_1)$ and one of the branches of $W^s(p_2)$. That is, $\mathcal O(x)$ is  
a {\it saddle connection}
from $p_1$ to $p_2$.

A sequence of saddle equilibrium points, $p_1, p_2, \cdots, p_n$, $p_{n+1}=p_1$ is
said to form a {\it heteroclinic cycle} if there is a saddle connection $\Gamma_i$ from
$p_i$ to $p_{i+1}$, $i=1,2, \cdots, n$. Heteroclinic cycles are prevalent in population 
dynamics \cite{hofbauer_evolutionary_1998, afraimovich_robust_2010}  and occur naturally in low codimension
bifurcation theory \cite{guckenheimer2013nonlinear}.
Fig. 1(a) shows an example of a 3-cycle.  Orbits starting from initial conditions sufficiently 
near $p_1$ in the domain bounded by $\Gamma := \cup \Gamma_i$ will 
shadow the piecewise smooth curve $\Gamma$  for some time. 
Additionally, if the heteroclinic cycle is stable, then $\varphi_t(x)$ will be attracted to $\Gamma$
as $t \to \infty$, in much the same way that nearby orbits are attracted to limit cycles ---
except that as $t$ increases, the orbit spends a larger and larger fraction of time near 
the three fixed points. 

A sufficient condition for the stability of a heteroclinic cycle connecting saddle fixed points
$p_1, \cdots, p_n$ is
\begin{align}
    \frac{- \lambda_s^{(i)}}{\lambda_u^{(i)}}>1 \qquad \mbox{for all } i
    \label{eq:stable_condish}
\end{align}
where $\lambda_s^{(i)}$ is the eigenvalue in the stable direction of $Df(p_i)$ and 
$\lambda_u^{(i)}$ is the eigenvalue in the unstable direction. Condition (\ref{eq:stable_condish}) 
implies that compression is greater than expansion at each $p_i$.
This condition can in fact be relaxed to requiring only that
the product of the ratios $- \lambda_s^{(i)}/\lambda_u^{(i)}$ be $>1$. See e.g. Refs.~\cite{krupa_asymptotic_1995, krupa_robust_1997, melbourne_heteroclinic_1989}.

\begin{figure}
\centering
\includegraphics[width=0.65\linewidth]{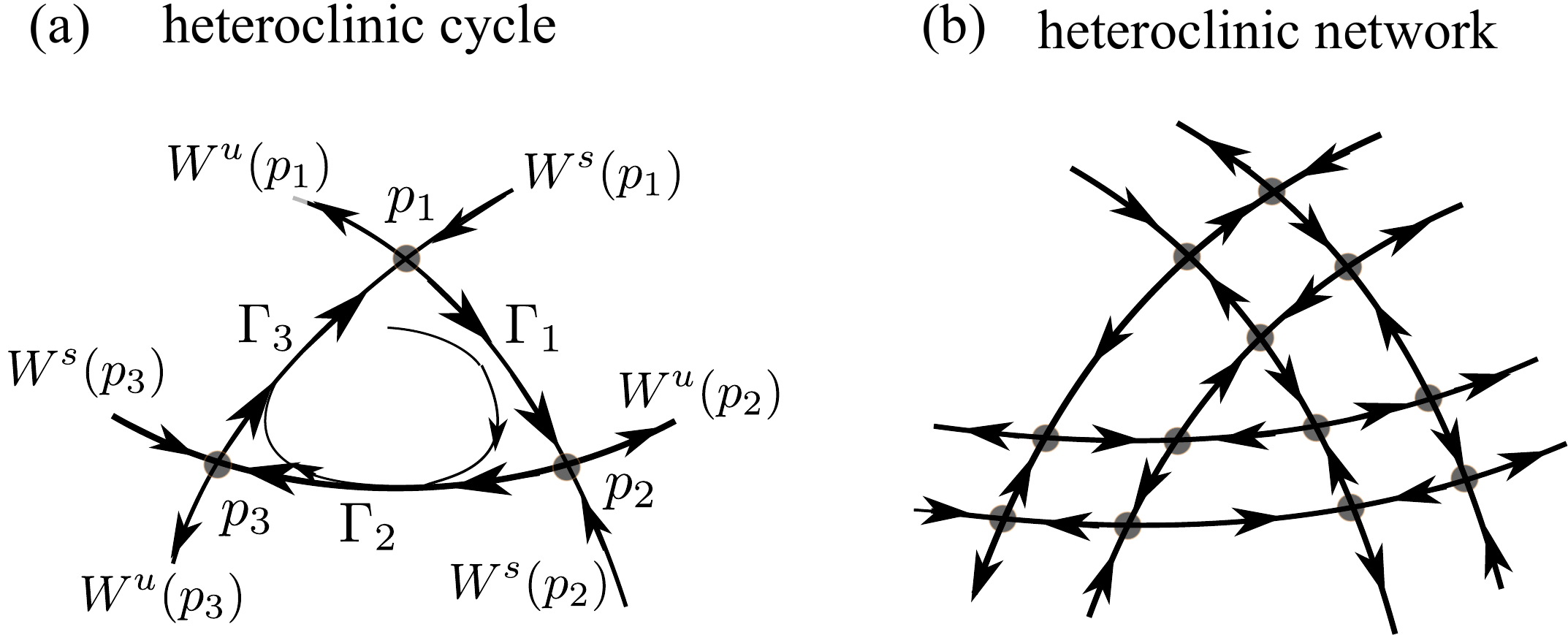}% Here is how to import EPS art
\caption{\label{fig:intro_fig2} (a) A heteroclinic cycle is created by connecting three heteroclinic orbits. (b) Many hyperbolic fixed points connected by heteroclinic orbits create a heteroclinic network.}
\end{figure}

More general than heteroclinic cycles is the idea of a {\it heteroclinic network}, loosely defined
to be the union of a finite collection of saddle fixed points joined by saddle connections.
An example of such a network is shown in Fig. 1(b). For further discussion and examples,
see e.g. Refs.~\cite{kirk_competition_1994, meiss2007differential, guckenheimer_structurally_1988, field_stationary_1991}. 

It is easy to see that in a heteroclinic network on $\mathbb R^2$,
 the phase space is divided into cells bounded by heteroclinic cycles. 
Imposing stability conditions on each cycle as above, the dynamical picture can be
described as follows: Each orbit is trapped in a cell; it is attracted to 
the heteroclinic cycle that bounds the cell, lingering at each one of the saddle fixed
points which act as metastable states before moving to the next. 
The set of metastable states visited depends on
 initial conditions, but they are visited in a cyclical order without exception.
Dwell times at each metastable state, hence the time to 
 complete each cycle, increase without bound as time goes to infinity.

Heteroclinic networks for 2D flows as described above are already models of switching 
behavior, but the dynamics they describe are very special. In particular, branching cannot
occur, and the dynamics are not chaotic. 

A number of authors have introduced random noise to heteroclinic networks of the type
above to induce branching behavior. Ref.~\cite{rabinovich_transient_2008}, for example,
constructs dynamical systems containing such networks to model orbits moving between 
metastable cognitive states in the brain. Ref.~\cite{bakhtin_neural_2012} used small noise perturbations
near saddle fixed points to model decision making. There is a fair amount of rigorous
theory on randomly perturbed heteroclinic networks \cite{armbruster_noisy_2003, bakhtin_noisy_2011, bakhtin_small_2010, bakhtin_rare_2022}. 
We postpone  further discussion of these results as our approach is orthogonal: 
The seemingly random switching behavior in our model comes not from
the use of random noise but from chaotic dynamics as we now describe.

\medskip \noindent
{\it Creation of branching behavior: qualitative picture}

\smallskip
Chaotic behavior in our model comes from the nontrivial intersection of stable and unstable
manifolds in heteroclinic networks. As discussed above, flows in 2D cannot support 
chaotic dynamics. The lowest phase dimension for which chaotic heteroclinic dynamics
can occur is 3D for continuous time and 2D for discrete time. For simplicity, we will work 
with the latter, though many of the ideas are not confined to two phase dimensions. 

Specifically, we start from a flow $\varphi_t$  on $\mathbb R^2$ with a heteroclinic network,  fix
a small number $t_0>0$, and let $F_0$ be the time-$t_0$-map $\varphi_{t_0}$. 
It is an easy fact that the saddle fixed points and their stable and unstable manifolds 
for $F_0$ are identical to those of $\varphi_t$. We describe below a procedure for
modifying $F_0$ along a saddle connection to obtain a smooth map 
$F: \mathbb R^2 \to \mathbb R^2$ with nontrivial intersection of
stable and unstable manifolds. The qualitative ideas behind this procedure 
are standard in the dynamical systems literature \cite{palis_hyperbolicity_1995}, 
but the {\it quantitative control} to follow are, to our knowledge, novel.

Let $p_1$ and $p_2$ be saddle fixed points. We assume there is 
a heteroclinic orbit $\Gamma_1$ of $\varphi_t$ connecting $p_1$ to $p_2$.
For simplicity, we assume, as in Fig. 2(a), that $p_1, p_2$ and $\Gamma_1$ lie 
on a horizontal line, while $W^s(p_1)$, $W^u(p_2)$ are vertical lines.
For the map $F_0$, a {\it fundamental domain} of $W^u(p_1)$ on $\Gamma_1$ 
is a segment $S \subset \Gamma_1$ with the property that the left endpoint of $S$
is mapped to the right endpoint under $F_0$, so that $F_0(S) \cap S = \emptyset$ 
modulo end points and
$\Gamma_1 = \cup_{n \in \mathbb Z} \ F_0^n(S)$.

\begin{figure}[!ht]
\centering
\includegraphics[width=0.6\linewidth]{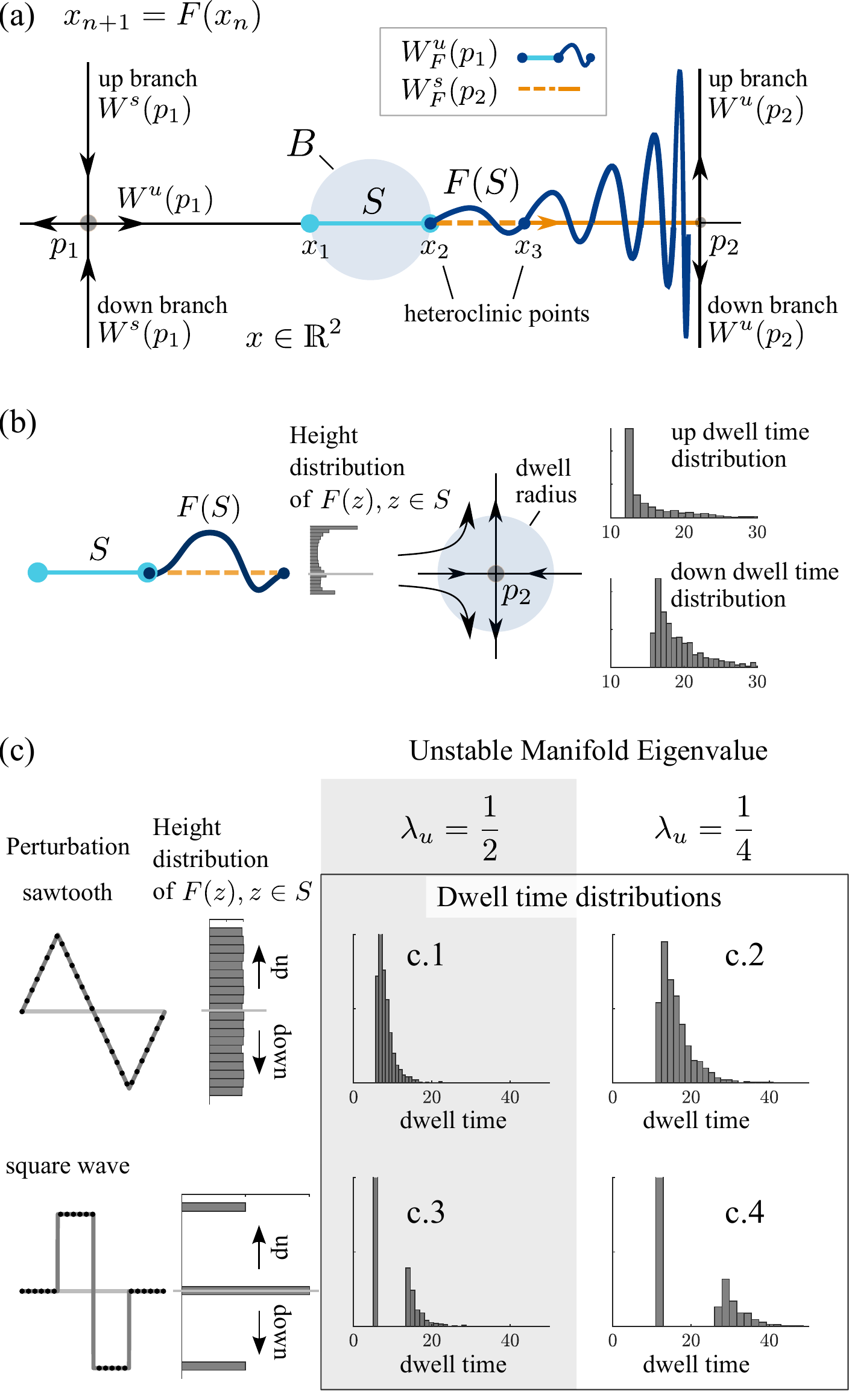}% Here is how to import EPS art
\caption{\label{fig:dwell_times} 
(a) A local perturbation near the segment $S$ is applied to the time-$t_0$-map $F_0$
of the 2D flow to break the saddle connection between fixed points $p_1$ and $p_2$,
creating nontrivial intersections of the stable and unstable manifolds, hence branching
behavior in the resulting map $F$.  Here $S=[x_1, x_2]$ (cyan)
is a fundamental domain in
$W^u(p_1)$; $F_0(S) = [x_2, x_3]$; $F(S)$ is the sinewave curve
joining $x_2$ to $x_3$ obtained by pushing points in $F_0(S)$ vertically up or down;
$F=F_0$ outside of the shaded region $B$.
(b) More detailed analysis of the perturbation in (a):
height distribution of points in $S$ one step later and statistics of subsequent dwell times 
near $p_2$, starting from uniform distribution on $S$.
 (c) Two examples of $g$ (sawtooth and square wave) are shown to
result in  different height distributions and dwell time statistics.
  A smaller unstable eigenvalue $\lambda_u$ at $p_2$, e.g.,
$\lambda_u = \frac14$, leads to longer dwell times  
(c.2, c.4) than $\lambda_u = \frac12$ (c.1, c.3).}
\end{figure}

Now fix a ball $B$ centered at the midpoint of $S$ with radius half the length of $S$, 
and let $F$ be such  that
\begin{alignat*}{2}
F & =F_0  \qquad && \text{outside of } B, \\
F & =F_0+g  \qquad && \text{on } B
\end{alignat*}
%
% \begin{eqnarray*}
% F=F_0 & \qquad & \mbox{outside of } B, \\
% F=F_0+g & \qquad & \mbox{on } B
% \end{eqnarray*}
%
where for all $x \in B$, $F_0(x)+g(x)$ is assumed to lie directly above or below $F_0(x)$.
That is, we perturb $F_0$ into $F$ by pushing the images of points in $B$ up or down.
See Fig. 2(a), where $F(S)$ has the shape of a biased sinewave. Since $F=F_0$ outside of $B$,
 the $F^n$-images of $F(S)$  for $n=1,2,3,\cdots$ are stretched vertically,  compressed
horizontally and pressed against $W^u(p_2)$ as $n$ increases.

Let $W_F^u(p_1)$ denote the unstable manifold of $F$ through $p_1$ 
(to distinguish it from the corresponding object for $F_0$), and let $W_F^s(p_2)$
denote the stable manifold of $F$ through $p_2$.
 Using shorthands such as
$S=[x_1, x_2]$ to refer to the horizontal segment from $x_1$ to $x_2$,
and letting $F(x_2)=F_0(x_2)=x_3$ (see Fig. 2(a)), it is easy to check that
\begin{eqnarray*}
\mbox{right branch of } W_F^u(p_1) & = & [p_1, x_1] \ \bigcup \ \cup_{n \ge 0} F^n(S)\ , \\
\mbox{left branch of } W_F^s(p_2) & = & [x_2, p_2] \ \bigcup \ \cup_{n \ge 0} F^{-n}([x_2, x_3]) \ .
\end{eqnarray*}
Thus, if $g$ is nontrivial, meaning if $F$ maps some points in $S$ above $\Gamma_1$
and some points below, then $W_F^u(p_1)$ and $W_F^s(p_2)$ intersect nontrivially
and branching behavior is created: Points in $B$ that are mapped to  locations above
$\Gamma_1$ will eventually follow the up-branch of $W_F^u(p_2)$, while those 
mapped to locations below $\Gamma_1$ will follow the down-branch of $W_F^u(p_2)$.

The chaotic behavior or sensitive dependence on initial condition that ensues
stems from the fact that for a point $x$ near $p_1$, it is hard to predict 
whether its trajectory under $F$ will follow the up or down-branch of $W_F^u(p_2)$. 
For $x \in \Gamma_1$ strictly to the right of $p_1$, $F^n(x)$ will be in $S$ for some  $n>0$.
In the scenario depicted in Fig. 2(a), if it lies on the left two-thirds of 
$S$, it will follow the up-branch, and if it lies on the right third, it will go down.
But even though its future is predetermined, the precise location of $F^n(x)$ is hard
to predict when $x$ is very closer to $p_1$ and $n$ is large. 
 For $x \not \in \Gamma_1$ near $p_1$ and on the right side of $W^s(p_1)$, 
due to the vertical compression and horizontal expansion near $p_1$, there will be an 
$n>0$ such that $F^n(x)$ is very closer to $S$. 
The same reasoning as above
then applies, with the precise location of $x$ 
playing a role in the ``decision" at $p_2$.

\medskip \noindent
{\it Examples of quantitative control}

\smallskip
Fig. 2(b) elucidates the branching properties that result from the perturbation
described in Fig. 2(a). Fixing a disk around $p_2$, we define {\it dwell time} near $p_2$
to be the duration of time a trajectory spends in the disk. Observe that points that are farther 
from $\Gamma_1$ will have shorter dwell times than those closer to $\Gamma_1$.
In this example, the distribution of dwell times for trajectories that eventually follow
the up- and down-branches are different due to the up-down bias in the map $g$.
A histogram showing the heights of $\{F(z_i)\}$ for a large collection of points 
$\{z_i\}$ evenly spaced on $S$ is also shown.
The higher probability to be above $\Gamma_1$ (hence to follow the up-branch
of $W^u_F(p_2)$ is evident, as is the fact that among those $F(z_i)$ above $\Gamma_1$,
the center of mass is farther from $\Gamma_1$, resulting in shorter dwell times.

Consider a line segment $S'$ slightly above $S$: It is easy to see that 
the height distribution of $F(z_i)$ and dwell times (not shown) will be even more biased
than starting from $S$. For $S'$ higher still, all trajectories will follow the up-branch.
A similar analysis can be carried out for all points near $S$.

Two other examples of $g$ are shown in Fig. 2(c). The top and bottom
panels show transition statistics for trajectories perturbed with a sawtooth perturbation
versus a square wave perturbation. (One can approximate these functions 
with smooth maps without affecting too much the distributions of interest.)
The sawtooth distribution results in the heights of 
$F(z_i)$ being uniformly distributed, whereas the square wave distribution pushes some points
far from $\Gamma_1$ and leaves others close, resulting in distinct clusters. We illustrate
also, in the right panels (c.1-c.4), the effect of the unstable eigenvalues on dwell times.
These plots confirm that the shapes of dwell time distributions are determined by 
the height distributions of $F(z_i)$ as well as the eigenvalues at the saddle
fixed point: weaker expansion at $p_2$ pushes points away  
more slowly, thereby increasing dwell times.

%%%%%%%%%%%%%%%%%%%%%%%%%%%%%%%%%%
%%%%%%%%%%%%%%%%%%%%%%%%%%%%%%%%%%%
%%%%%%%%%%%%%%%%%%%%%%%%%%%%%%%%%%%%
\section{\label{sec:section3} Heteroclinic network built to reproduce four-state  \textit{C. elegans} 
behavioral data}

\textit{C. elegans}, nematodes that grow to about $1$mm in length, are a well-studied model organism in neuroscience \cite{white_structure_1986, chalasani_dissecting_2007, kimata_thermotaxis_2012, varshney_structural_2011, kato_global_2015}. Whole-brain calcium imaging techniques record the activity of most \textit{C. elegans}' neurons while real or fictive behavioral states are simultaneously observed \cite{nguyen_whole-brain_2016, kato_global_2015, randi_measuring_2020}.
Because their movements are mediated by waves of contractions of bands of muscles that run
the full length of the body, 
\textit{C. elegans} exhibit only a limited set of behaviors, predominantly forward crawling, turns, reversals, and quiescence \cite{kato_global_2015, nichols_global_2017, linderman_hierarchical_2019}. 
The low-dimensional neural dynamics of \textit{C. elegans} cluster into different regions in PCA space that correspond to distinct behaviors \cite{nichols_global_2017, kato_global_2015}. 
Experimentalists have fit Markov models to the behavioral sequences observed,
treating transitions between behavioral states as well as the dwell times in each state as stochastic
\cite{nichols_global_2017, linderman_hierarchical_2019}.
Variations of behavior patterns  in different experimental conditions have also been documented \cite{linderman_hierarchical_2019, nichols_global_2017}.

In the last two sections of this paper we build dynamical systems 
in the form of chaotic heteroclinic networks 
(as described in Section II) to simulate \textit{C. elegans} data. 
Results reported from two different labs will be used to challenge the model. 
In this section we focus on results from Ref.~\cite{nichols_global_2017}.

\begin{figure}[ht!]
\centering
\includegraphics[width=0.6\linewidth]{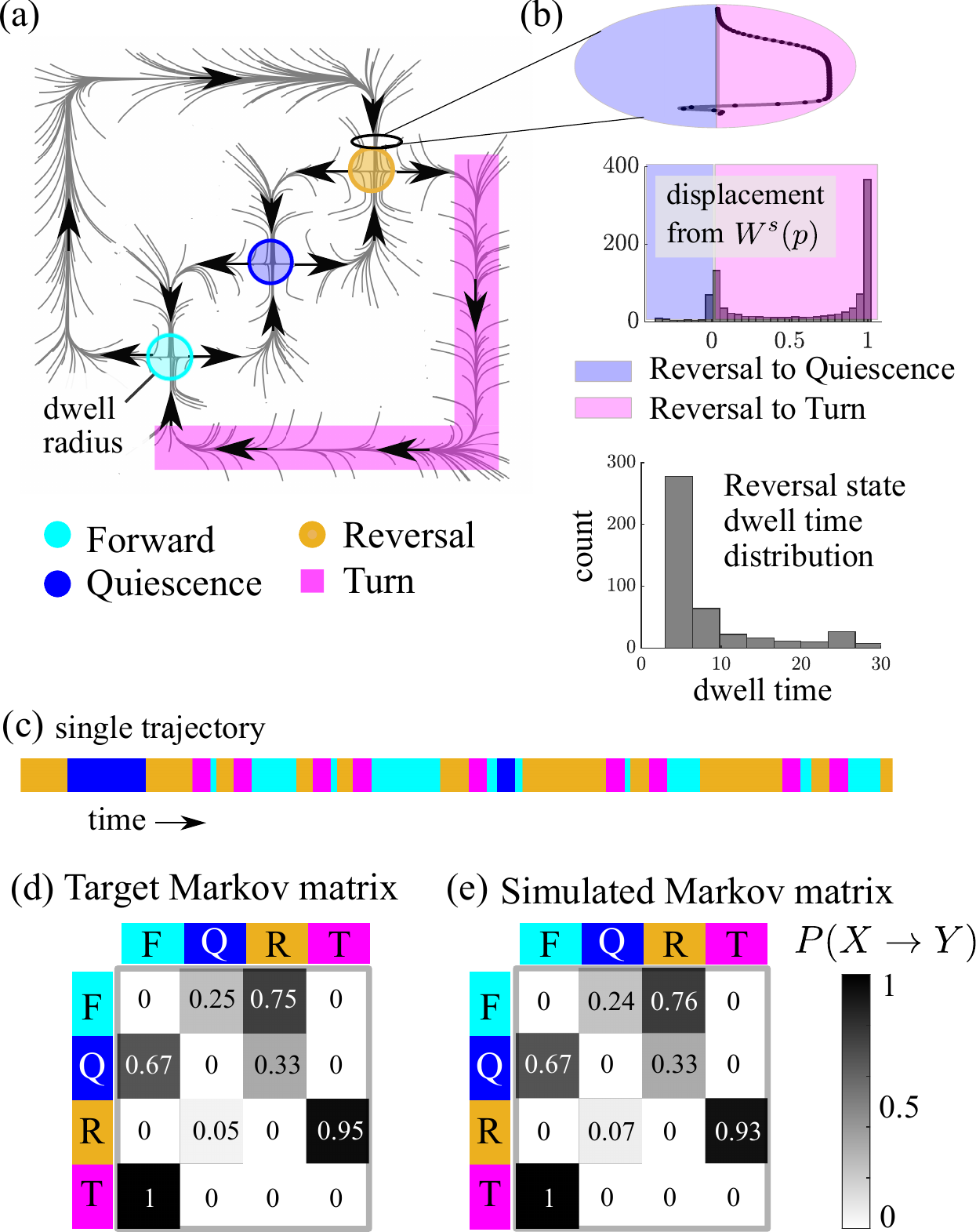}% Here is how to import EPS art
\caption{\label{fig:Nichols2017} (a) Phase portrait of a dynamical system containing a heteroclinic network intended to recreate \textit{C. elegans} transition statistics in Ref.~\cite{nichols_global_2017}. Unperturbed heteroclinic network flow lines are shown in grey.  Regions around fixed points correspond to \textit{C. elegans} behavioral states Forward (cyan), Quiescence (blue), and Reversal (yellow). The Turn (pink) behavioral state is treated as 
a transitional state and represented by a heteroclinic orbit. (b) Branching behavior is created
by perturbing the heteroclinic orbits prior to the fixed points. The shape of the perturbation determines left/right biases and dwell times. Reversal state dwell times is measured as the time spent within a certain radius of the reversal fixed point. (c) Example trajectory of behavioral states
(the short time intervals between states are removed from the time series).
(d) Target transition probabilities estimated from Ref.~\cite{nichols_global_2017}. (e) Probabilities from system simulation fit to match target probabilities in (d).}
\end{figure}

\bigskip
Ref.~\cite{nichols_global_2017} considered primarily four types of behaviors, abbreviated
as {\it forward, reversal, turn}, and {\it quiescence}. Treating these four behaviors as discrete states, 
we construct a chaotic heteroclinic
network (Fig.~\ref{fig:Nichols2017}(a)) containing three saddle fixed points that correspond 
to the forward, quiescence, and reversal states, and model the ``turn" state by a heteroclinic orbit 
connecting the reversal to the forward fixed point: turns in the data seem 
always to be preceded by reversal and followed by forward movements.  
As discussed in Section II, we first build a (nonchaotic) heteroclinc network 
 (in the form of a 2D flow) with 
the properties above. The chaotic network is obtained by applying 
perturbations to the flow-map along relevant 
heteroclinic orbits. Construction details are provided 
in the {\bf Appendix}.
Dwell times are defined to be time spent in a specified region around each one of the
saddle fixed points or a portion of the heteroclinic orbit in the case of turns.

For illustration, we pick one of the modifications of the flow-map,
namely the one illustrated in Fig.~\ref{fig:Nichols2017}(b), and elaborate on
the considerations that went into our choice of the perturbation function, $g$,
in the region shown. 
Experimental data show that following reversal, the turn state occurs
with very high probability ($0.95$ according to the Markov matrix in (d)). 
Assuming that the quiescence state is not dominant (which is often the case), 
most orbits arriving at reversal are from the upper branch
of the stable manifold $W^s(p)$ where $p$ is the fixed point corresponding to reversal. 
These two facts together suggest that $g$, the function that determines
how the upper branch of  $W^s(p)$ meets the unstable manifold from the fixed point
corresponding to the forward state, should have a strong bias. Accordingly,
we have chosen a perturbation that pushes most of the points to the right 
of $W^s(p)$ and eventually into the turn state.
 Ref.~\cite{nichols_global_2017} did not include data on how to
constrain the shape of $g$, giving only bias probabilities. In general, the shape
of  $g$ can be deduced from measured dwell time distributions. The $g$ shown
in Fig.~\ref{fig:Nichols2017}(b), e.g., produces
a bimodal distribution of displacements (or deviations from $W^s(p)$) one step later; 
see Section II for more detail.  Such a distribution (Fig.~\ref{fig:Nichols2017}(b)) implies 
that a large fraction of the trajectories is pushed quite far from $W^s(p)$  
and will have shorter dwell time at $p$,
while a smaller fraction stays very close to it and will have longer dwell times.

To show how a typical trajectory in the constructed dynamical system may look,  we
create a behavioral state time series by labeling the system at each point in time with either
one of the four identified states or as ``in-between", referring to times not in any of the four
sets. Transitions between states are
assumed to be relatively fast, making the ``in-between" category negligible. 
In the example shown in Fig.~\ref{fig:Nichols2017}(c), we have cut out the ``in-between" 
stretches showing only the transition dynamics between the four states.

A target Markov matrix containing transition statistics estimated from experimental data published in Ref.~\cite{nichols_global_2017} is shown in Fig.~\ref{fig:Nichols2017}(d).
These statistics were used to guide our model construction.
The corresponding statistics collected from our model 
are shown in Fig.~\ref{fig:Nichols2017}(e). 
Similarity in the matrix entries shows that perturbation functions can be tailored to 
produce simulated data that closely match the target statistics.

\medskip 
We digress here to point out that unlike Markov models, the dynamical models 
we propose have the capability to encode memory of past events. The following are two
examples to illustrate this flexibility:
(1) From the discussion above, we see that 
most orbits from forward to reversal leave reversal on the outside of the heteroclinic orbit connecting reversal to forward. If we set the strength of the contraction during the time
in the turn state to be weak, most orbits will remain on the outside of the loop at forward, 
making them more likely to switch back to reversal after a short visit to forward. 
A strong contraction along the heteroclinic orbit will pinch all orbits close
to the stable manifold of the fixed point at forward, resulting in less correlation between
past and future behaviors. (2) We have prescribed a function, $g$, 
for how the unstable manifold from forward meets the upper branch of $W^s(p)$. 
There is a corresponding function, $\hat g$, to connect 
the unstable manifold from quiescence to the {\it lower} branch of  $W^s(p)$.
By choosing $g$ and $\hat g$ to be different, one can cause dwell times and
 transition probabilities at reversal to depend on the previous step, i.e., whether
 the orbit came from forward or quiescence. Similar constructions can lead to 
 dependencies on the steps before.

\medskip
Another set of results from Ref.~\cite{nichols_global_2017} that we would like to replicate
concerns \textit{C. elegans}' transition statistics and dwell times under different experimental conditions: {\it prelethargus}
vs {\it lethargus}, and $10\%$ vs $21\%$ oxygen. The first are developmental stages:
It is known that during lethargus periods, \textit{C. elegans} are more likely
 to transition to quiescent behavior and to remain in a quiescent state for longer stretches of time with brief periods of activity \cite{nichols_global_2017}. Nematodes in the prelethargus stage, in contrast, are much more active.
Oxygen levels in the environment are also known to affect activity level:
\textit{C. elegans} in an atmospheric 21\% oxygen environment
are more likely to remain active 
than \textit{C. elegans} in a 10\% oxygen environment, which is amenable to sleep. 
This condition occurs naturally when worms socially aggregate, creating a preferred low oxygen environment
\cite{nichols_global_2017}.
Ref.~\cite{nichols_global_2017} contains detailed information on 
 dwell times
at and transition probabilities from the forward state under different
developmental and environmental conditions.

\begin{figure}[ht!]
\centering
\includegraphics[width=0.6\linewidth]{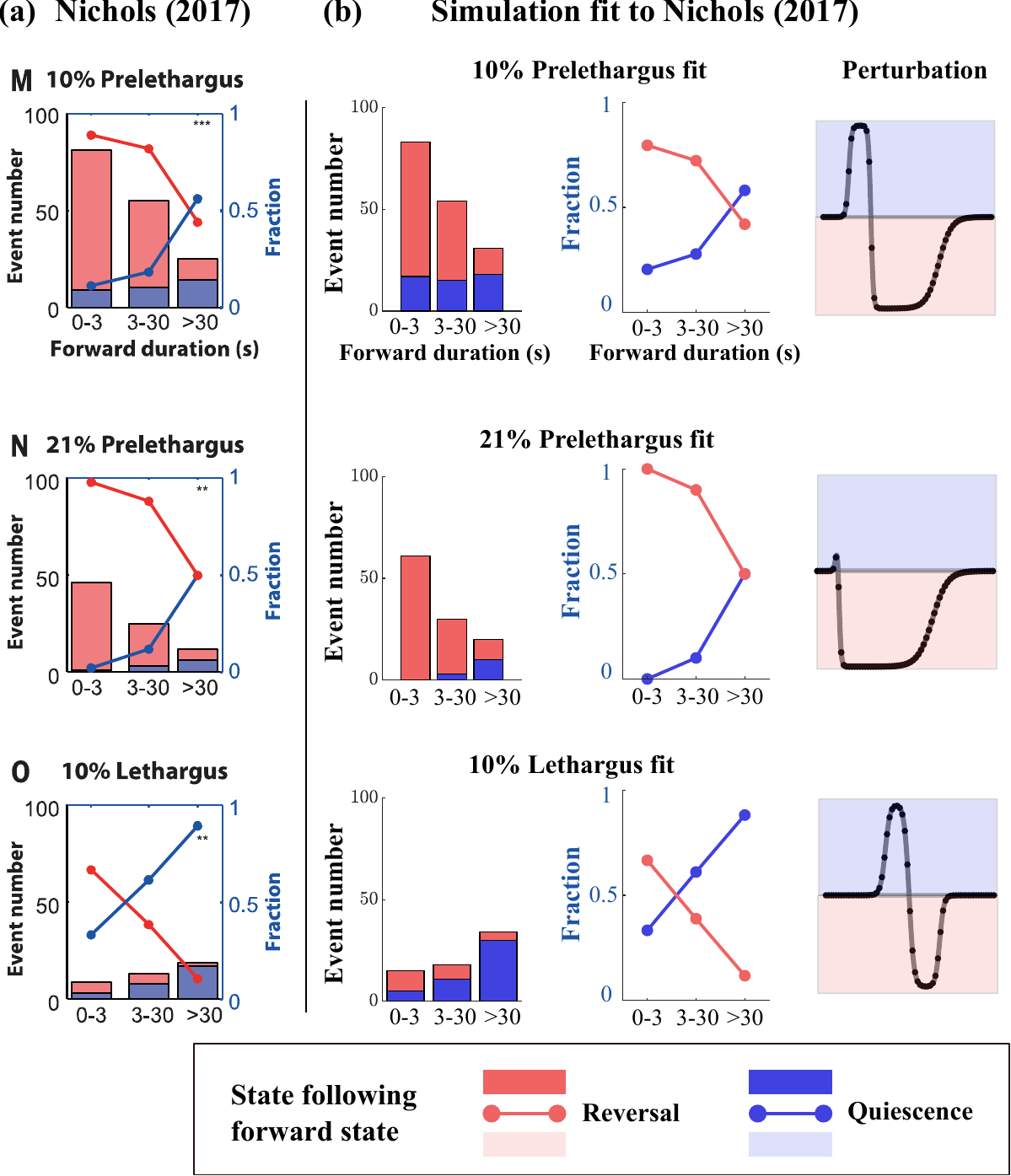}% Here is how to import EPS art
\caption{\label{fig:Nichols2017_2} (a) Dwell time and bias statistics for three different experimental conditions - (M) 10\% Prelethargus, (N) 21\% Prelethargus, and (O) 10\% Lethargus. (From Ref.~\cite{nichols_global_2017}. Reprinted with permission from AAAS) (b) Corresponding statistics produced by heteroclinic network shown in Fig.~\ref{fig:Nichols2017}(a). The first column is a histogram of transitions out of forward state after increasingly large dwell times. The second column is the fraction of transitions to the reversal versus quiescent state following forward crawling (middle two columns). Perturbation shapes used to generate these results is shown in the right most column.}
\end{figure}

Fig.~\ref{fig:Nichols2017_2}(a) reprints, with permission from AAAS, data
from Ref.~\cite{nichols_global_2017}. Each panel reports on data recorded
under a different set of conditions.  The heights of the bars in each histogram 
(ignoring colors for now) show the numbers of transitions out of the forward state within 
$0$-$3$, $3$-$30$ and $>30$ seconds after entering this state; they can be seen 
as a probability distribution of dwell times. 
Each bar in the histogram is further decomposed into a red and
a blue portion representing the number of transitions to reversal (red) or quiescence (blue). 
The red and blue lines superimposed on the histogram indicate the fractions of 
transitions of these two types. In the middle two columns, we present the corresponding
statistics collected from our heteroclinic network.

In the top row we fit transition statistics for \textit{C. elegans} in a prelethargus state with a $10\%$ oxygen environment. In this developmental stage \textit{C. elegans} are more active.
It is observed experimentally that 
under these conditions, about half of
the transitions out of forward occur in the first $3$ seconds. Moreover, when dwell time is $< 3$ sec., 
the \textit{C. elegans} are much more likely to switch from forward to reversal;
the longer the dwell time in the forward state, the more likely the \textit{C. elegans} will switch to quiescence. Overall, transitions from forward to reversal are much more prominent than transitions
into quiescence, 
though there is a nonnegligible fraction of the latter due to the low oxygen condition.
These data are reproduced in our dynamical systems model using the perturbation 
function $g$ shown in the rightmost panel. Comparing the middle two panels to the left, 
we see the strong quantitative resemblance of  model outputs to experimental data.

Next we compare the top two rows, both of which are for the prelethargus state but at different
oxygen levels. Transition statistics are qualitatively similar,
except that at $21\%$ oxygen,  transition probabilities to quiescence are a little bit
lower as high oxygen levels promote activity in \textit{C. elegans}. Our model is able to capture these relatively subtle distinctions using 
the perturbations shown.

Differences between the prelethargus and lethargus data (at 10\% oxygen) 
are far more substantial. In contrast to the top row, dwell times in the forward state
in the lethargus case (bottom row) are much longer,
nearly half of them more than 30 seconds. As before, the longer the time to transition, 
the more likely it will go into quiescence. This translates into a very significant difference
between the prelethargus and lethargus data: In the lethargus case, a good majority
of the transitions are to quiescence. 
Fig.~\ref{fig:Nichols2017_2} shows excellent agreement between data and model
in both cases.
 
We finish by observing that even though only four types of behaviors are considered,
a $4$-state Markov chain is not sufficient for capturing 
the type of results in Fig.~\ref{fig:Nichols2017_2}; our dynamical model with its
tunable parameters offers a simple way to recreate these and other behavioral characteristics.

\section{\label{sec:section4}Heteroclinic network built to reproduce eight-state \textit{C. elegans} behavioral data}

\begin{figure}[ht!]
\centering
\includegraphics[width=0.6\linewidth]{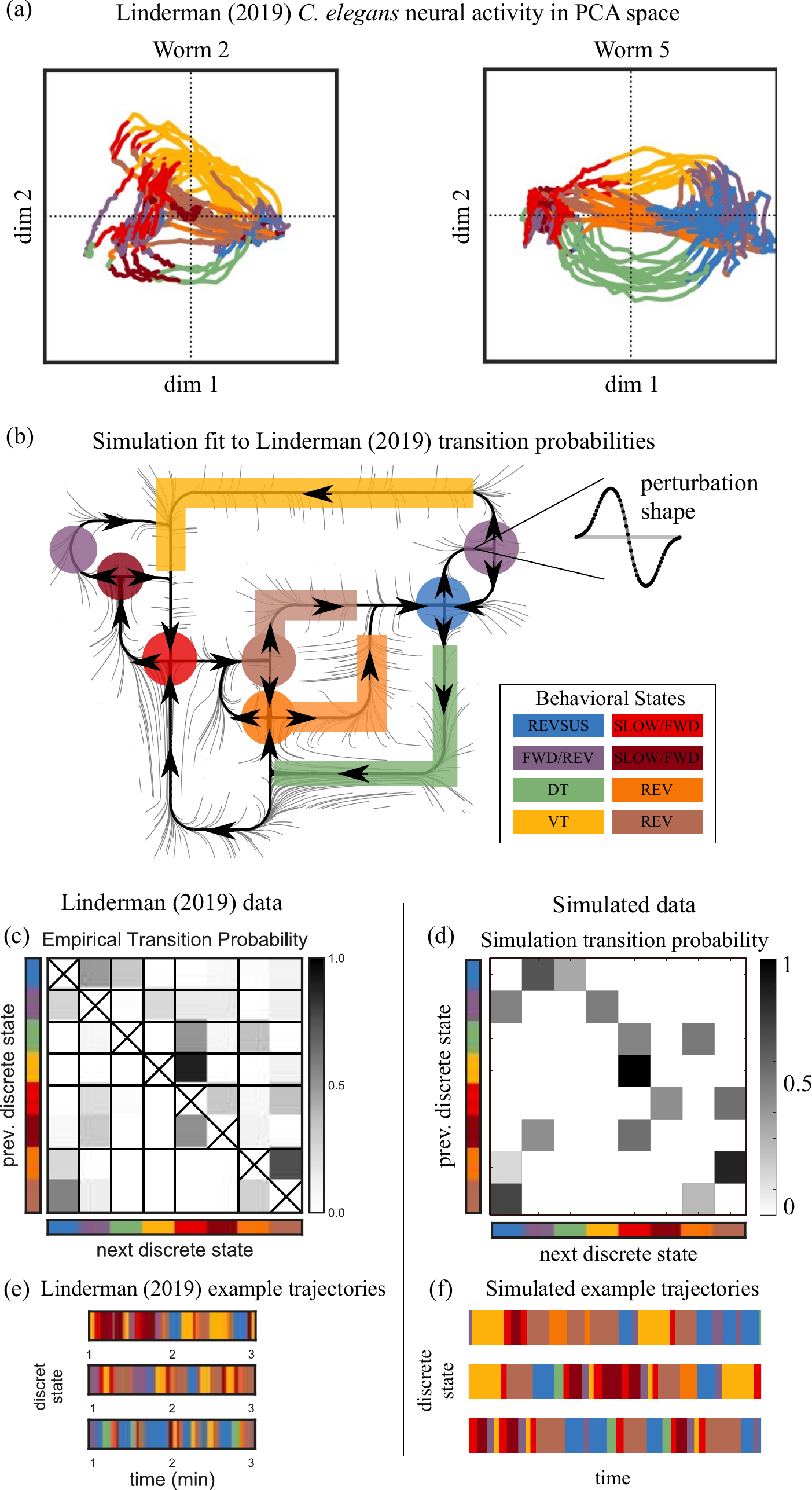}% Here is how to import EPS art
\caption{\label{fig:Linderman} (a) 2D projections in PCA space of trajectories
segmented into states shown in different colors (same color key as in (b)) for two different worms
(Reprinted from Ref.~\cite{linderman_hierarchical_2019} with permission). (b) Nonlinear dynamical system embedding a heteroclinic network with corresponding behavioral states. (c)
Markov matrix of behavioral state transition statistics reprinted from Ref.~\cite{linderman_hierarchical_2019}. (d) Transition probabilities from dynamical system simulated data. Markov matrix transition probabilities are fit to approximate Markov matrix in (c) by adjusting perturbation functions. (e) Example behavioral trajectories from Ref.~\cite{linderman_hierarchical_2019}. (f) Example behavioral trajectories simulated from nonlinear system.}
\end{figure}

%%%%%%%%%%%%%%%%%%%%%%%%%%%%%%%%
%%% Complex C. elegans %%%%%%%%%
%%%%%%%%%%%%%%%%%%%%%%%%%%%%%%%%
In other studies \cite{linderman_hierarchical_2019, kato_global_2015}, 
researchers have made more refined classifications of \textit{C. elegans} behavioral states than those used to fit the model in the previous section. 
Fig.~\ref{fig:Linderman}(a) is reproduced from Ref.~\cite{linderman_hierarchical_2019}
with permission. It shows a 2D projection (from high dimensional data) of trajectories 
segmented into states with characteristic dynamics shown in different colors for two different worms. 
Here the authors distinguished between 8 discrete behavioral states, among them 
{\it sustained reversals} (blue), {\it dorsal} and {\it ventral turns}
(green and yellow), {\it forward crawling} (red and crimson) and the {\it transition from forward back to
reversal} (orange and brown). As part of their study, they proposed an 8-state Markov model, 
the transition probabilities of which are shown in Fig.~\ref{fig:Linderman}(c) (reproduced with permission). 
Our aim in this section is to demonstrate that a chaotic heteroclinic network model can be built
to reproduce the more complex picture of \textit{C. elegans} behavior reported 
in this paper.

As in the last section, we start by building a 2D flow, representing each of the eight states 
by either a saddle fixed point or a heteroclinic orbit. Connections between nodes are based on  Figure~\ref{fig:Linderman}(a),(c); a diagram is shown in Figure~\ref{fig:Linderman}(b).
The color code is as in Figure~\ref{fig:Linderman}(a) and (c); one of the states (purple) 
appears twice as trajectory segments in that state are found
 in two distinct sets of circumstances. To keep the dynamics simple, we have opted to 
 capture only  the more significant transitions, omitting several that occur with very low probability.
Modifications are then made on the time-$t$-map of the flow to create nontrivial
intersections of stable and unstable manifolds resulting in a chaotic heteroclinic
network. The perturbation functions used are chosen with the aim of 
reproducing the transition probabilities in Figure~\ref{fig:Linderman}(c). 

In addition to transition probabilities, Ref.~\cite{linderman_hierarchical_2019}, 
Fig. 6(b), provides some guidance on the duration of time spent in each state.
Though numerical values are not provided, and it is impossible to deduce exact dwell times
from these figures, they do offer clues on dwell times. 
We have incorporated these clues into 
our dynamical model through the choice of unstable eigenvalues 
at saddle fixed points and by adjusting the speeds 
of travel along heteroclinic orbits.
By specifying different eigenvalues and perturbation functions for each fixed point, as well as specifying different speeds along the corridors connecting fixed points, we are able to tune the model to replicate both transition statistics and dwell times in the 8-state \textit{C. elegans} behavioral data.

Model outputs of the resulting dynamical system are collected, and the likelihood of its trajectories 
visiting state $j$ immediately after state $i$ is tabulated and shown in Fig.~\ref{fig:Linderman}(d). Though not perfect, these dynamical transition probabilities --- generated by a purely deterministic
dynamical system --- show good resemblance to the stochastic transition probabilities in Fig.~\ref{fig:Linderman}(c). Finally, we present in Fig.~\ref{fig:Linderman}(f) three example trajectories from our dynamical network, in analogy with the three trajectories from  
 Ref.~\cite{linderman_hierarchical_2019} (Fig.~\ref{fig:Linderman}(e)).
The variable itineraries and patterns of transitions in the three trajectories shown 
exemplify the possibilities and unpredictable outcomes typical in chaotic dynamical systems.

%%%%%%%%%%%%%%%%%%%%%%%%
%%% Discussion %%%%%%%%%
%%%%%%%%%%%%%%%%%%%%%%%%
\section{\label{sec:discussion}Discussion}

{\it Stochastic switching in biology.} 
The paradigm of having a finite number of dominant states and switching between
them in a seemingly random way in the absence of external stimuli is relevant
beyond \textit{C. elegans}. It
has been observed in biological and neural activity in many organisms,
and  Markov models have been used to capture these behaviors.
For example, mice have been found to have different behavioral 
modules, the switching dynamics between which can be modeled with 
a hidden Markov model \cite{wiltschko_mapping_2015}. 
Larval zebrafish have been observed to transition between discrete behavior states 
following a Markov model \cite{sharma_point_2018}.
Fruit fly behavior is known to evolve along a low-dimensional attractor with periodic
and nonperiodic behavioral sequences \cite{berman_mapping_2014}. 
Other transitions are likely to be caused or facilitated by biological events, but 
when a mechanistic understanding is out of reach due to high complexities in 
the underlying biology, many authors have, as a first attempt, idealized such transitions as random.
For example, activity in the ventromedial prefrontal cortex (vmPFC) of humans appears
to transition stochastically through a series of discrete states that correspond to affective 
experiences; these state transitions have been modeled with an HMM \cite{chang_endogenous_2021}. The neural activity in the gustatory cortex appears to
transition between discrete states stochastically, a fact some researchers have
attributed to noise-induced variability \cite{miller_stochastic_2010}.
On the phenomenological level, these and many other examples of switching dynamics 
have a great deal in common with those considered in this paper.

\medskip 
{\it Dynamical vs. stochastic models.} Seemingly random switching has been 
observed in dynamical models. For example, models of spiking neurons are known to support 
multiple semi-stable states visited by the system in ways that resemble random transitions \cite{mazzucato_dynamics_2015}.  It is also known that random switching in (purely deterministic)
dynamical systems can be generated by the  addition of random noise. Ref.~\cite{rabinovich_transient_2008}, for example, posits that transient cognitive dynamics result from metastable cognitive states which, when linked together, produce stable heteroclinic channels. 
Sequential decision-making can be modeled with stable heteroclinic channels with noise 
added to the system. See also Refs.~\cite{bakhtin_small_2010, bakhtin_neural_2012, bakhtin_noisy_2011}. 
Differential dwell times near fixed points and random switching between heteroclinic cycles 
generated by random noise have been studied by a number of authors; see e.g. Refs.~\cite{stone_random_1990, armbruster_noisy_2003}. 

As an alternative to Markov models or dynamical models perturbed by random noise, we propose
in this paper the use of chaotic heteroclinic networks to model seemingly random switching 
behavior. In our heteroclinic network models, dynamical complexity alone --- without the need for 
stochastic manipulations --- produces the variability and unpredictable behavior observed.
One argument in favor of dynamical
models could be that they are more natural, in the sense that physical and
biological phenomena are 
not really stochastic in nature. One could also argue that 
dynamical systems models are more flexible: In Sections III and IV, we have demonstrated
that chaotic heteroclinic networks can be constructed to emulate  \textit{C. elegans} 
data in two influential papers. Our models accurately reproduce not only basic
transition properties between semi-stable states, such as dwell times and branching biases 
(Figs.~\ref{fig:Nichols2017}, \ref{fig:Nichols2017_2} and \ref{fig:Linderman}), but more subtle behaviors, including exit distributions 
and dependence of transition probabilities on dwell time (Fig.~\ref{fig:Nichols2017_2}). We have also
explained how one could, through network design, influence the degree to which
memory of prior events will be retained. 

This is not to suggest that stochastic models cannot
possess similar capabilities; they can, but the models will have to be somewhat more
complicated than $n$-state Markov chains where $n$ is the number of behavioral states. 
We have shown that dynamical models that are relatively easy to construct 
are capable of recreating a range of behaviors,
with excellent quantitative control  achieved through  tunable parameters. 
Hence we propose them as a viable alternative to purely stochastic models. 

\medskip 
{\it Remarks on global network construction.} 
As discussed, we start by constructing a 2D flow which we perturb to obtain the
desired switching behavior. Among the many ways to construct such a flow,
we have found the following to be easy to tune and to analyze: First we
lay down all the saddle fixed points in $\mathbb R^2$, then we connect them
 with stable and unstable manifolds for the 2D flow. Each stable (unstable) manifold
 is assumed to comprise a vertical and a horizontal segment with corners rounded
 (as in Fig.~5(b)).
 In a neighborhood of this system of ``train tracks" of stable and unstable manifolds, 
 we build the vector field to be contractive to trap all nearby orbits.
Where stable and unstable manifolds cross in $\mathbb R^2$, ``overpasses" can
be introduced. This geometry does not  impose constraints 
on the architecture of the network.

\medskip 
{\it Extension to a larger class of dynamical models.}
For clarity of exposition, we have limited ourselves in this paper to heteroclinic networks 
that are perturbations of time-$t$-maps of 2D flows. There are many ways to enlarge
this class of models. An obvious extension is to allow higher dimensions,
so that the stable and unstable manifolds of the saddle fixed points can have dimension
greater than one.
In the rest of the Discussion, however, we would like to focus on a different extension,
to a broad class of dynamical systems we call
 {\it generalized heterclinic networks}:

Continuing to consider switching behavior among a finite collection $\mathcal S$ of identifiable, 
semi-stable states, an extension that we think is natural --- both
from the perspective of modeling and from the dynamical systems theory point of view 
--- is to consider internal dynamics within each of the states in $\mathcal S$. 
Specfically, the saddle fixed points in classical heteroclinic networks can be  replaced 
by hyperbolic invariant sets such as hyperbolic periodic orbits, normally hyperbolic invariant tori supporting quasi-periodic behavior, or Smale's horseshoes \cite{smale_differentiable_1967}.  
In the case of a periodic cycle, for example, trajectories approaching the cycle will 
be attracted to it, stay with it for some time, and then veer away after
going around the cycle a seemingly random number of times. Stable and unstable manifolds of hyperbolic invariant sets are well defined \cite{palis_hyperbolicity_1995}
and their intersections will define the switching dynamics.

\subsection*{Acknowledgments}

This work was supported by the National Science Foundation (LSY, grant no. 1901009) as well as the National Science Foundation Mathematical Sciences Postdoctoral Research Fellowship (MM, award no. 2103239).
 Part of this work was done when LSY held a visiting position at the Institute for Advanced Study.   

 We would like to thank Annika Nichols and Scott Linderman for permission to reproduce figures in their respective publications, Refs.~\cite{nichols_global_2017} and \cite{linderman_hierarchical_2019}.  We would also like to thank Charles Fieseler for his valuable comments on our manuscript.

\subsection*{Data Availability Statement}
    The data that support the findings of this study are available within Refs.~\cite{nichols_global_2017} and \cite{linderman_hierarchical_2019}.

% %%%%%%%%%%%%%%%%%%%%%%%%
% %%%  Appendix  %%%%%%%%%
% %%%%%%%%%%%%%%%%%%%%%%%%

\appendix

\section*{Appendix: Heteroclinic network construction}\label{appx:A}

We provide here more analytical details on the heteroclinic network construction used
in this paper. The methods outlined below are generic and can serve as templates for network constructions
elsewhere.

\begin{figure}[ht!]
\centering
\includegraphics[width=0.6\linewidth]{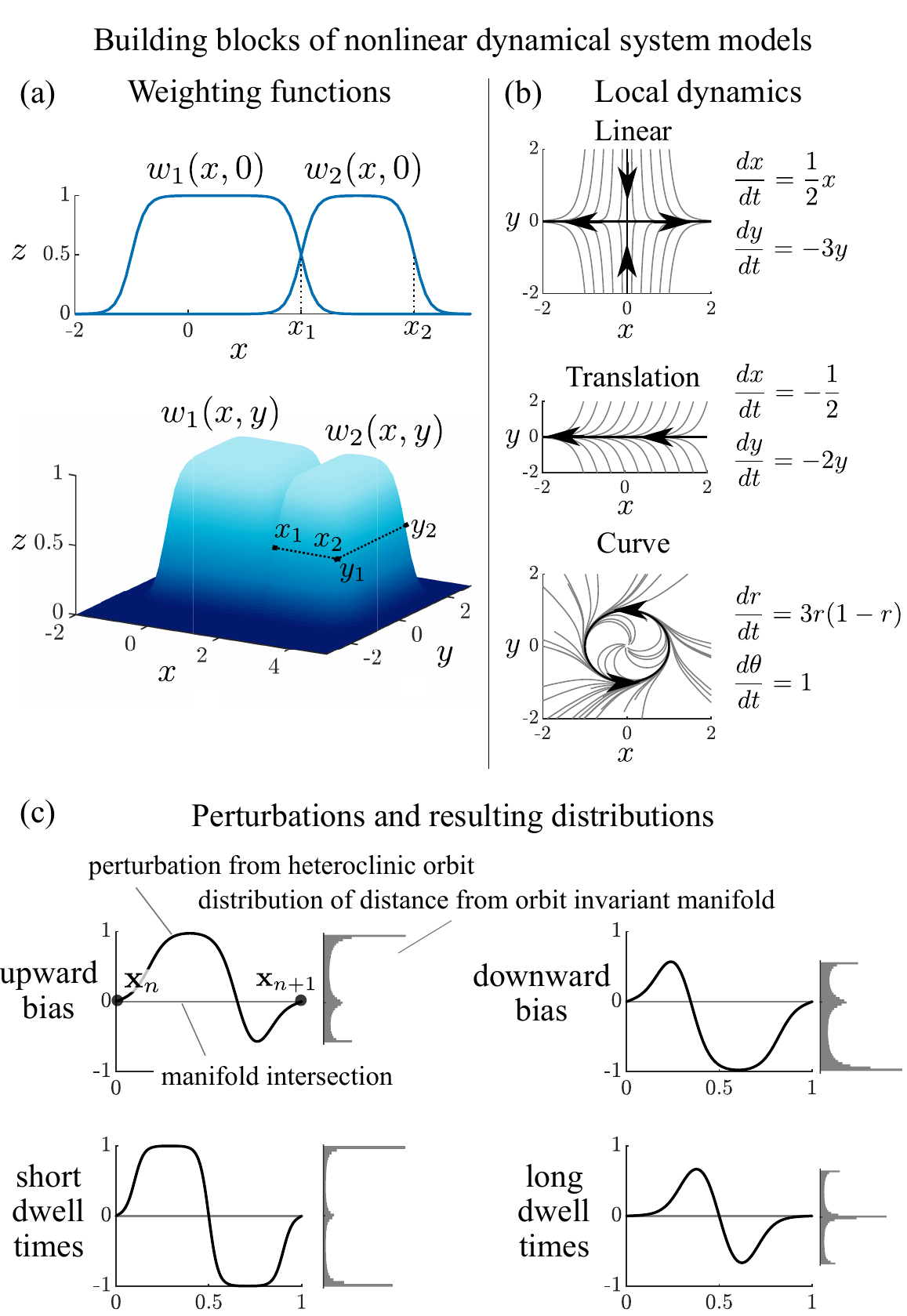}
\caption{\label{fig:building_blocks} (a) Weighting functions for local dynamics. (b) Local dynamics used to create heteroclinic network. (c) Types of perturbations. Points perturbed above the manifold intersection travel upward while points perturbed below travel downward. Points perturbed far from the intersection have short dwell times near the subsequent fixed point while points perturbed infinitesimally from the intersection have long dwell times near the subsequent fixed point.}
\end{figure}

We build a nonlinear dynamical system embedding a heteroclinic network by setting different functions to be dominant in different regions of phase space.
By interpolating local dynamics we create a global nonlinear landscape. The global dynamics are
\begin{align}\label{eq:global_dyn}
    \frac{d \b{x}}{dt} &= \sum_{i=1}^m w_i(\b{x}) f_i(\b{x}).
\end{align}
where the local dynamics $f_i(\b{x})$ are weighted by $w_i(\b{x})$,
\begin{align}
\begin{split}
   w(x,y) = \frac{1}{4}(\tanh[s(x-x_1)] - \tanh[s(x-x_2)])*\\
    (\tanh[s(y-y_1)] - \tanh[s(y-y_2)]).
\end{split}
\end{align}
The weighting functions, $w_i(\b{x})$, weight the dynamics, $f_i(\b{x})$, highly in each local dynamics' region of influence, $(x,y) \in [x_1,x_2] \times [y_1,y_2]$ (Fig.~\ref{fig:building_blocks}(a)). The sharpness of the transition between the dynamics in different regions is controlled by hyperbolic tangent slope $s$. In the limit as $s \rightarrow \infty$ the system becomes a  piecewise function.

We use three types of local functions $f_i(\b{x})$ to construct the heteroclinic network: linear dynamics (Eq.~\ref{eq:linear_dyn}), transversal dynamics (Eq.~\ref{eq:trans_dyn}), and rotational dynamics (Eq.~\ref{eq:rot_dyn}). Stitched together, these three types of dynamics can create flexible heteroclinic networks.

The structure of the $f_i(\b{x})$ producing linear dynamics is
\begin{align}
\begin{split}\label{eq:linear_dyn}
    \frac{dx}{dt} &= \lambda_1 x\\
    \frac{dy}{dt} &= \lambda_1 y.
\end{split}
\end{align}
The $x$ and $y$ axes form the beginning of the stable and unstable manifolds emanating from the fixed point at the origin. The signs of $\lambda_1$ and $\lambda_2$ control which axis is the stable versus unstable manifold; the magnitude of  $\lambda_1$ and $\lambda_2$ determine how quickly the system moves toward and away from the fixed point. The location of the fixed point can be shifted to any location $p_i \in \mathds{R}^2$.

The structure of the transversal dynamics is
\begin{align}
\begin{split}\label{eq:trans_dyn}
        \frac{dx}{dt} &= c\\
    \frac{dy}{dt} &= -a (y-b)
\end{split}
\end{align}
where constants $a$, $b$, and $c$ determine the strength of the attracting manifold at $y=b$ and the direction and speed of movement along it. Transversal dynamics can also be constructed to move vertically.

The dynamics for the rotations are easier to represent in polar coordinates and can be expressed as
\begin{align}
\begin{split}\label{eq:rot_dyn}
    \frac{dr}{dt} &= a r(b - r)\\
    \frac{d\theta}{dt} &= c
\end{split}
\end{align}
where constants $a$, $b$, and $c$ determine the strength of the attracting radius at $r=b$ and the rotational speed and direction. Polar coordinates can be transformed into Cartesian coordinates and the location of the radius center can be relocated to any $p_i \in \mathds{R}^2$.

Figure~\ref{fig:building_blocks}(b) shows the three types of local dynamics which are weighted by weights $w_i(\b{x})$ in order to link the stable and unstable manifolds emanating from fixed points.
These building blocks give us great flexibility in creating heteroclinic orbits which comprise the heteroclinic network.

After a continuous time version of the heteroclinic network is constructed (Eq.~\ref{eq:global_dyn}), a discrete mapping in constructed from the continuous dynamics
\begin{align}
    \b{x}_{n+1} = \b{x}_n + dt \sum_{i=1}^m w_i(\b{x}_n) f_i(\b{x}_n).
\end{align}
Perturbations are subsequently applied to the heteroclinic orbits prior to each stable fixed point in this discrete time dynamical system.  Perturbation functions perturb the trajectory off of the heteroclinic orbit and take the form
\begin{align}
g(x) =& 
\begin{cases} 
      g_+(x) & x \geq 0 \\
      g_-(x) &  x<0
   \end{cases}\\
   \begin{split}
   g_+(x) =& \frac{1}{2}(\tanh[s_1(x-L_1)]+\tanh[s_1(x+L_1)])\\
   & - \tanh[s_2(x-(\frac{1}{2}+b))]\\
   & \frac{1}{2}(\tanh[s_3(x-(1-L_2))]-\tanh[s_3(x-(1+L_2))])\end{split}\\
   g_-(x) &= g_+(-x)
\end{align}
where parameters $L_1$, $L_2$, $s_1$, $s_2$, $s_3$, and $b$ control the distribution of points that are perturbed above versus below the heteroclinic orbit. These parameters also control the distance from the orbit that these points are displaced.
 Figure~\ref{fig:building_blocks}(c) shows perturbations with different parameter values that result in more points perturbed above versus below the heteroclinic orbit and close to versus far from the heteroclinic orbit. The width of the perturbation function is scaled to span one timestep of the discrete map and is applied only once preceeding each fixed point. The amplitude of the perturbation is also scaled and affects dwell times. Without perturbations, trajectories would converge to heteroclinic orbits in the network and follow stable manifolds to fixed points, spending longer amounts of time at each subsequent fixed point as trajectories converge closer to the invariant manifolds comprising the heteroclinic orbits. With the use of the perturbation functions, in addition to tuning the parameters for the linear, transversal, and rotational dynamics, dwell time distributions and transition probabilities can be imposed onto the network.

%\nocite{*}
%\bibliographystyle{abbrvnat}
%\bibliography{Lai-Sang,Nathan, Zhuo-Cheng,me, Nonlinear_control_Koopman}% Produces the bibliography via BibTeX.

\end{document}